\documentclass[11pt,bezier]{article}
\usepackage{amsmath}
\usepackage{amsfonts,amsthm,amssymb}
\usepackage{amsfonts}
\usepackage{graphicx}
\usepackage{epstopdf}
\usepackage{caption}
\usepackage{amssymb}
\newtheorem{lem}{Lemma}[section]
\newtheorem{thm}[lem]{Theorem}
\newtheorem{cor}[lem]{Corollary}

\theoremstyle{definition}

\begin{document}
\title{Super connected direct product of graphs and cycles\footnote{The research is supported by National Natural Science Foundation of China (11861066).}}
\author{Jiaqiong Yin, Yingzhi Tian\footnote{Corresponding author. E-mail: YinJiaQiong6@163.com (J. Yin); tianyzhxj@163.com (Y. Tian).} \\
{\small College of Mathematics and System Sciences, Xinjiang
University, Urumqi, Xinjiang, 830046, PR China}}

\date{}
\maketitle
\begin{sloppypar}

\noindent{\bf Abstract }
The topology of an interconnection network can be modeled by a graph $G=(V(G),E(G))$. The connectivity of graph $G$ is a parameter to measure the reliability of corresponding network. Direct product is  one important graph product. This paper mainly focuses on the
super connectedness of direct product of graphs and cycles.

The connectivity of $G$, denoted by $\kappa(G)$, is the size of a minimum vertex set $S\subseteq V(G)$ such that $G-S$ is not connected or has only one vertex. The graph $G$ is said to be super connected, simply super-$\kappa$, if every minimum vertex cut is the neighborhood of a vertex with minimum degree. The direct product of two graphs $G$ and $H$, denoted by $G\times H$, is the graph with vertex set $V(G \times H) = V (G)\times V (H)$ and edge set $E(G \times H) = \{(u_{1} ,v_{1} )(u_{2} ,v_{2} )|\ u_{1}u_{2} \in E(G), v_{1}v_{2} \in E(H)\}$. In this paper, we give some sufficient conditions for direct product $G\times C_{n}$ to be super connected, where $C_{n}$ is the cycle on $n$ vertices. Furthermore, those sufficient conditions are best possible.


\section{Introduction}

For a  simple graph $G$ with vertex set $V(G)$ and edge set $E(G)$, $u, v\in V(G)$ are adjacent if $uv\in E(G)$.  The set of all vertices adjacent to $u$ is called the $neighborhood$ of $u$ in $G$, denoted by $N_{G}(u)$. The $degree$ of $u$, denoted by $d_{G}(v)$, is $|N_{G}(u)|$. The $minimum$ $degree$ of $G$ is $\delta(G)=min\{d_{G}(v)| v\in V(G)\}$. For a vertex set $S\subseteq V(G)$, if $G-S$ is not connected, then $S$ is a vertex cut of $G$.  We known that only complete graphs do not have vertex cuts.  If $G$ is not a complete graph, then the $connectivity$ of $G$, denoted by $\kappa(G)$, is the size of a minimum vertex cut of $G$. Otherwise, $\kappa(G)=|V(G)|-1$.
The $edge$ $connectivity$ of $G$, denoted by $\kappa'(G)$, is the the size of a minimum  edge set  $F\subseteq G$  such that $G-F$ is not connected.  $K_n$, $K_{m,n}$ and $C_n$ are used to denote complete graph, complete bipartite graph and cycle, respectively. We follow Bondy and Murty \cite{Bondy} for undefined notation and terminology.

The topology of an interconnection network can be modeled by a graph $G=(V(G),E(G))$,
where $V(G)$ represents the set of nodes and $E(G)$ represents the set of communication links in the network. The connectivity $\kappa(G)$ of $G$  can be used to measure the reliability and fault tolerance of the network. In generally, the larger $\kappa(G)$ is, the more reliable the network is. It is well known that $\kappa(G)\leq\delta(G)$. The graph $G$ with $\kappa(G)=\delta(G)$ is called $maximally$ $connected$, simply $max$-$\kappa$. For the maximally connected graphs, it is believed that the graphs with the smallest number of minimum vertex cuts are more reliable than the others. Boesch in \cite{Boesch} proposed the concept of super connected graph. If every minimum vertex cut is a neighborhood of some vertex of $G$ with minimum degree, then the graph $G$ is said to be $super$ $connected$, simply $super$-$\kappa$. By the definitions, a super connected graph is also maximally connected. The converse is not always true. For example, $C_n\ (n\geq6)$  is maximally connected but not super connected.

The direct product of two graphs $G$ and $H$, denoted by $G\times H$, is the graph with vertex set $V(G \times H) = V (G)\times V (H)$ and edge set $E(G \times H) =\{(u_{1} ,v_{1} )(u_{2} ,v_{2} )|\ u_{1}u_{2} \in E(G), v_{1}v_{2} \in E(H)\}$. Weichsel \cite{Weichsel} proved that the direct product $G\times H$ of two nontrivial connected graphs $G$ and $H$ is connected if and only if at least one of  $G$ and $H$ are not  bipartite.

We list some  results on the edge connectivity of direct product of graphs as follows. Some bounds on the edge connectivity of the direct product of graphs were given by \v{S}pacapan  in \cite{Bresar}. Cao,  Brglez,  \v{S}pacapan and  Vumar \cite{Cao} determined the edge connectivity of direct product of a nontrivial graph and a complete graph. In \cite{Spacapan}, \v{S}pacapan not only obtained the edge connectivity of direct product of two  graphs, but also characterized the structure of a minimum edge cut in the direct product of two graphs.

This paragraph will list some results on the connectivity of direct product of graphs. Some bounds on the connectivity of the direct product of graphs were also given by \v{S}pacapan in \cite{Bresar}.  Mamut and Vumar \cite{Mamut} proved that the connectivity of the direct product of two complete graphs $K_m $ and $K_n$ is $(m-1)(n-1)$, where $m\geq n\geq2$. In \cite{Guji}, Guji and Vumar proved that the connectivity of the direct product of a complete graph $K_n\ (n\geq3)$ and  a bipartite graph $G$  is $min\{n\kappa(G),(n-1)\delta(G)\}$,  and furthermore, the authors also conjectured that this is true for all nontrivial graph $G$. Later, Wang and Wu \cite{Wang&Wu} and Wang and Xue \cite{Wang&Xue} independently confirmed this conjecture. Wang and Yan \cite{Wang&Yan} determined the connectivity of $G\times K_2$. Recently, Sonawane and Borse \cite{Sonawane} determined the connectivity of the direct product of graphs and cycles.

The  results on the super connected direct product graphs are presented in the following.
Guo, Qin and Guo \cite{Guo} proved that for a maximally connected bipartite graph $G$, $G\times K_n$ ($n\geq 3$) is super connected. In \cite{Wang&Shan}, the authors generalized this result by showing that for a maximally connected  nonbipartite graph $G$, $G\times K_n$ ($n\geq 3$) is super connected. In \cite{Zhou}, Zhou completely characterized the super connected direct product of a nontrivial graphs $G$ and a complete graph $K_{n}$ ($n\geq3$), that is,   $G\times K_{n}$ is not super connected if and only if either $\kappa(G\times K_{n})=n\kappa(G)$ or $\kappa(G\times K_{n})\cong K_{m,m}\times K_{3}\ (m\geq1)$.  Wu and Tian \cite{Wu} studied the super connected direct product of paths, cycles and cycles.

Motivated by the results above, especially those in \cite{Sonawane}, we will study super connected direct product of graphs and cycles. In the next section, we will present a key lemma, which will be used in the proof of our main results in Section 3. Conclusion will be given in the last section.

\section{A key lemma}

%


In \cite{Sonawane}, Sonawane and Borse constructed a graph $\widetilde{G}_{n}$ from a connected bipartite graph $G$ as follows. Let $(X,Y )$ be a bipartition of $G$ and $n\geq 2$ be an integer and $X_{i}=\{(x,i):x\in X\}$ and $Y_{i}=\{(y,i):y\in Y\}$ for $i=1,2,...,n$. Let $H_{i}$ and $H_{i}^{'}$ be graphs isomorphic to $G$ with bipartitions $(X_{i},Y_{i})$ and $(X_{i+1},Y_{i})$, respectively for each $i\in\{1,2,...,n\}$. Let
\[
\begin{array}{ll}
\widetilde{G}_{n}=\bigcup_{i= 1}^n {\left({{H_i\cup H_{i}^{'}}}\right)}.
\end{array}
\]

\begin{figure}[htbp]
  \centering
  \includegraphics[width=13cm]{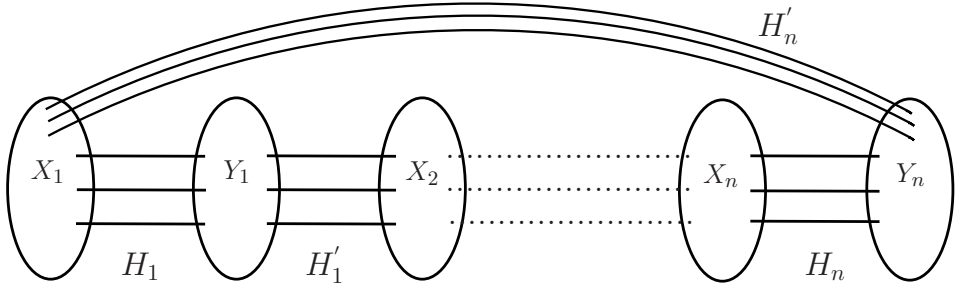}\\
  \caption{The graph $\widetilde{G}_{n}$}
  \label{1}
\end{figure}

Sonawane and Borse \cite{Sonawane} determined the connectivity of $\widetilde{G}_{n}$ as follows.

\begin{thm} (\cite{Sonawane}) $\kappa(\widetilde{G}_{n})=min\{n\kappa(G), 2\delta(G)\}$, where $n\geq 2$ is an integer.
\end{thm}

In the following, we will give a key lemma in this paper, which will be used to prove the main results in the next section. Furthermore, the sufficient condition for $\widetilde{G}_{n}$ to be super connected in this lemma is best possible.

\begin{lem} $\mathrm{(}${\textrm{\bf A key lemma}}$\mathrm{)}$
Let $G$ be a connected bipartite graph with  bipartition $X$ and $Y$, and let $n\geq 3$ be an integer. If  $|X|\geq\delta(G)+1$, $|Y|\geq\delta(G)+1$  and $\kappa(G)>\frac{2}{n}\delta(G)$, then $\widetilde{G}_{n}$ is super-$\kappa$.
\end{lem}

\noindent{\bf Proof.} By Theorem 2.1 and $\kappa(G)>\frac{2}{n}\delta(G)$, we have $\kappa(\widetilde{G}_{n})=2\delta(G)$. By contradiction, assume that $\widetilde{G}_{n}$ is not super-$\kappa$, then there is a vertex cut $S$ with $|S|=2\delta(G)$ such that $\widetilde{G}_{n}-S$ is not connected and has no isolated vertices. Let $D_{1},D_{2},...,D_{r}\ (r\geq2)$ be the  components of $\widetilde{G}_{n}-S$. Then $|D_{i}|\geq2$ holds for each $i\in\{1,...,r\}$. Denote $S_{X_{i}}=S\cap X_{i}$, $S_{Y_{i}}=S\cap Y_{i}$, $X_{i}^{'}=X_{i}-S_{X_{i}}$ and $Y_{i}^{'}=Y_{i}-S_{Y_{i}}$ for $i=1,2,...,n$.

Since each vertex $x_i$ in $X_i$ has at least $\delta(G)$ neighbors in both $Y_{i-1}$ and $Y_{i}$,  each vertex $y_j$ in $Y_j$ has at least $\delta(G)$ neighbors in both $X_{j}$ and $X_{j+1}$, we have the following claim.

\noindent{\bf Claim 1.} If $D_{k}\cap X_{i}^{'}\neq\emptyset$, $|S_{Y_{i-1}}|<\delta(G)$ and $|S_{Y_{i}}|<\delta(G)$, then both $D_{k}\cap Y_{i-1}^{'}\neq\emptyset$ and $D_{k}\cap Y_{i}^{'}\neq\emptyset$ hold. Similarly, if $D_{k}\cap Y_{j}^{'}\neq\emptyset$, $|S_{X_{j}}|<\delta(G)$ and $|S_{X_{j+1}}|<\delta(G)$, then both $D_{k}\cap X_{j}^{'}\neq\emptyset$ and $D_{k}\cap X_{j+1}^{'}\neq\emptyset$ hold.

In the following, we consider two cases.

\noindent{\bf Case 1.}  $|S_{X_{i}}|<\delta(G)$ and $|S_{Y_{j}}|<\delta(G)$ for any $i,j\in\{1,2,...,n\}$.

By Claim 1, we have $D_{k}\cap X_{i}^{'}\neq\emptyset$ and $D_{k}\cap Y_{j}^{'}\neq\emptyset$ for  $k\in\{1,...,r\}$ and $i,j\in\{1,...,n\}$.
Since $|S|=(|S_{X_{1}}|+|S_{Y_{1}}|)+(|S_{X_{2}}|+|S_{Y_{2}}|)+\cdot\cdot\cdot+
(|S_{X_{n}}|+|S_{Y_{n}}|)<n\kappa(G)$, we have $|S_{X_{i}}|+|S_{Y_{i}}|<\kappa(G)$ for some $i\in\{1,2,...,n\}$. Note that $H_i=(X_i, Y_i)$ is isomorphic to $G$. Then $\widetilde{G}_{n}[X_i'\cup Y_i']=H_i-S_{X_{i}}\cup S_{Y_{i}}$ is connected. Thus
$\widetilde{G}_{n}-S$ is connected, contradicting to the assumption.

\noindent{\bf Case 2.} There is an $i\in\{1,2,...,n\}$ such that $|S_{X_{i}}|\geq\delta(G)$ or there is a $j\in\{1,2,...,n\}$ such that $|S_{Y_{j}}|\geq\delta(G)$.

Assume, without loss of generality, that $|S_{X_{1}}|\geq\delta(G)$.

\noindent{\bf Subcase 2.1.} $|S_{X_{i}}|<\delta(G)$ for any $i\in\{2,3,...,n\}$ and $|S_{Y_{j}}|<\delta(G)$ for any $j\in\{1,2,...,n\}$.

If $D_k\cap X_{1}'\neq\emptyset$ for some $k\in\{1,...,r\}$ , then by Claim 1,  $D_{k}\cap X_{i}^{'}\neq\emptyset$ and $D_{k}\cap Y_{j}^{'}\neq\emptyset$ for all $i,j\in\{1,...,n\}$. Otherwise, $D_{k}\cap X_{i}^{'}\neq\emptyset$ and $D_{k}\cap Y_{j}^{'}\neq\emptyset$ for all  $i\in\{2,...,n\}$ and $j\in\{1,...,n\}$.


Since $|S_{Y_{1}}|+(|S_{X_{2}}|+|S_{Y_{2}}|)+\cdot\cdot\cdot+(|S_{X_{n}}|+|S_{Y_{n}}|)
<n\kappa(G)-|S_{X_{1}}|\leq(n-1)\kappa(G)$, we have $|S_{X_{i}}|+|S_{Y_{i}}|<\kappa(G)$ for some $i\in\{2,3,...,n\}$.  Then $\widetilde{G}_{n}[X_i'\cup Y_i']=H_i-S_{X_{i}}\cup S_{Y_{i}}$ is connected. Thus $\widetilde{G}_{n}-S$ is connected, which contradicts to  the assumption.

\noindent{\bf Subcase 2.2.} There is an $i'\in\{2,3,...,n\}$ such that $|S_{X_{i'}}|\geq\delta(G)$.

Since $|S|=2\delta(G)$, we have  $|S_{X_{1}}|=\delta(G)$ and $|S_{X_{i'}}|=\delta(G)$.

If $i'=2$ or $n$, by symmetry, assume $i'=2$, then
$\widetilde{G}_{n}[Y_{2}'\cup X_{3}'\cup\cdot\cdot\cdot\cup X_{n}'\cup Y_{n}']=\widetilde{G}_{n}[Y_{2}\cup X_{3}\cup\cdot\cdot\cdot\cup X_{n}\cup Y_{n}]$ is connected. Furthermore, $\widetilde{G}_{n}[X_{2}'\cup Y_{2}'\cup X_{3}'\cup\cdot\cdot\cdot\cup X_{n}'\cup Y_{n}'\cup X_{1}']$ is connected. Since $\widetilde{G}_{n}-S$ has no isolated vertices, each vertex in $Y_1$ has at least one neighbor in $X_1'$ or $X_2'$. Therefore, $\widetilde{G}_{n}-S$ is connected, a contradiction.

If $i'\neq2$ and $i'\neq n$, then by $\widetilde{G}_{n}[Y_{1}\cup X_{2}\cup\cdot\cdot\cdot\cup Y_{i'-1}]$ and $\widetilde{G}_{n}[Y_{i'}\cup X_{i+1}\cup\cdot\cdot\cdot\cup X_{n}\cup Y_{n}]$ are connected, we obtain that $\widetilde{G}_{n}-S$ is connected, contradicting to the assumption.

\noindent{\bf Subcase 2.3.} There is a $j'\in\{1,2,...,n\}$ such that $|S_{Y_{j'}}|\geq\delta(G)$.

Since $|S|=2\delta(G)$, we have $|S_{X_{1}}|=\delta(G)$ and $|S_{Y_{j'}}|=\delta(G)$.

If $j'=1$ or $n$, by symmetry, assume $j'=1$, then
$\widetilde{G}_{n}[X_{2}'\cup Y_{2}'\cup\cdot\cdot\cdot\cup X_{n}'\cup Y_{n}']=\widetilde{G}_{n}[X_{2}\cup Y_{2}\cup\cdot\cdot\cdot\cup X_{n}\cup Y_{n}]$ is connected. Thus $\widetilde{G}_{n}-S=\widetilde{G}_{n}[Y_{1}'\cup X_{2}'\cup Y_{2}'\cup\cdot\cdot\cdot\cup X_{n}'\cup Y_{n}'\cup X_{1}']$ is connected, a contradiction.

If $j'\neq1$ and $j'\neq n$, then by $\widetilde{G}_{n}[Y_{1}\cup X_{2}\cup\cdot\cdot\cdot\cup X_{j'}]$ and $\widetilde{G}_{n}[X_{j'+1}\cup Y_{j'+1}\cup\cdot\cdot\cdot\cup X_{n}\cup Y_{n}]$ are connected, we obtain that $\widetilde{G}_{n}-S$ is connected, which contradicts to  the assumption.

Since all cases lead to contradiction, the proof is thus complete. $\Box$

\section{Main Results}

Motivated by the connectivity of the direct product of graphs and cycles, we will obtain some sufficient conditions for the direct product of graphs and cycles to be super connected. Considering four cases arising from whether $G$ is bipartite or not and $n$ is even or odd, Sonawane and Borse \cite{Sonawane} obtained the connectivity of the direct product of graphs and cycles in the following four theorems.

\begin{thm} (\cite{Sonawane}) Let $G$ be a connected bipartite graph and $n\geq3$ be an odd integer. Then $\kappa(G\times C_{n})=min\{n\kappa(G), 2\delta(G)\}$.
\end{thm}

\begin{thm} (\cite{Sonawane}) Let $G$ be a connected bipartite graph and  $n\geq4$ be an even integer. Then the graph $G\times C_{n}$ has two isomorphic components each with connectivity $min\{\frac{n}{2}\kappa(G), 2\delta(G)\}$.
\end{thm}

\begin{thm} (\cite{Sonawane}) Let $G$ be a connected non-bipartite graph and $n\geq4$ be an even integer. Then $\kappa(G\times C_{n})=min\{\frac{n}{2}\kappa(G\times K_{2}), 2\delta(G)\}$.
\end{thm}

\begin{thm} (\cite{Sonawane}) Let $G$ be a connected non-bipartite graph and $n\geq5$ be an odd integer. Then $min\{\frac{n-1}{2}\kappa(G\times K_{2}), 2\delta(G)\}\leq\kappa(G\times C_{n})\leq min\{\frac{n+1}{2}\kappa(G\times K_{2}), 2\delta(G)\}$.
\end{thm}

Similarly, we consider four cases to study the super connectedness of the direct product of graphs and cycles in the following. The cycle $C_{n}$ of length $n \ (\geq3)$ is denoted by $C_{n}=\langle1,2,...,n\rangle$.

\begin{thm}
Let $G$ be a connected bipartite graph with bipartition $X$ and $Y$, and let $n\geq3$ be an odd integer. If $|X|\geq\delta(G)+1$, $|Y|\geq\delta(G)+1$ and $\kappa(G)>\frac{2}{n}\delta(G)$, then $G\times C_{n}$ is super-$\kappa$.
\end{thm}

\noindent{\bf Proof.} Let $V_{i}=\{(v,i):v\in V(G)\}$ for $i=1,2,...,n$. Then $V(G\times C_{n})=\cup_{i= 1}^n {V_{i}}$. Let $X_{\frac{i+1}{2}}=\{(x,i):x\in X\}$ and $Y_{\frac{n+i}{2}}=\{(y,i):y\in Y\}$ for $i=1,3,...,n$. Let $X_{\frac{n+i+1}{2}}=\{(x,i):x\in X\}$ and $Y_{\frac{i}{2}}=\{(y,i):y\in Y\}$ for $i=2,4,...,n-1$. Then $V_{i}=X_{\frac{i+1}{2}}\cup Y_{\frac{n+i}{2}}$ for $i=1,3,...,n$ and $V_{i}=X_{\frac{n+i+1}{2}}\cup Y_{\frac{i}{2}}$ for $i=2,4,...,n-1$. Let
$H_{i}$ and $H_{i}^{'}$ be subgraphs of $G\times C_{n}$ induced by $X_{i}\cup Y_{i}$ and $X_{i+1}\cup Y_{i}$, respectively for $i=1,2,...,n$. Then $G\times C_{n}=\cup_{i= 1}^n {\left({{H_i\cup H_{i}^{'}}}\right)}$. The graph $G\times C_{n}$ is shown in Figure 2. For each $i$, $H_{i}$ and $H_{i}^{'}$ isomorphic to $G$. By Lemma 2.2, $G\times C_{n}$ is super-$\kappa$. $\Box$

\begin{figure}[htbp]
  \centering
  \includegraphics[width=13cm]{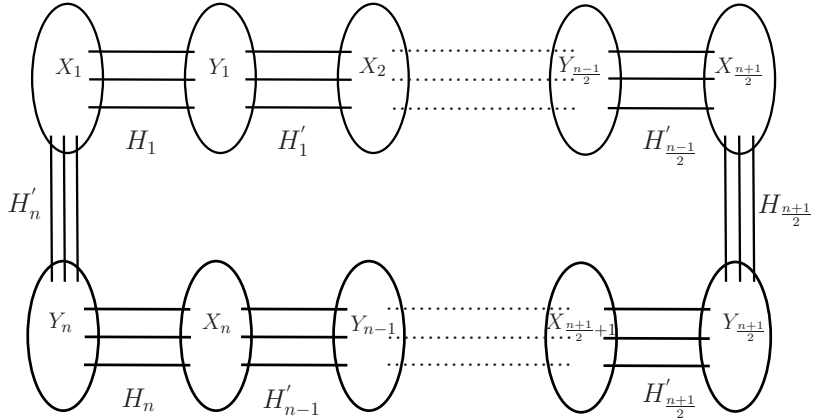}\\
  \caption{The graph $G\times C_{n}$ when $G$ is bipartite and $n$ is an odd}
  \label{1}
\end{figure}

\begin{thm}
Let $G$ be a connected bipartite graph with bipartition $X$ and $Y$, and let $n\geq6$ be an even integer. If $|X|\geq\delta(G)+1$, $|Y|\geq\delta(G)+1$ and $\kappa(G)>\frac{4}{n}\delta(G)$, then the two isomorphic components $G_{1}$ and $G_{2}$ of $G\times C_{n}$ are both super-$\kappa$.
\end{thm}

\noindent{\bf Proof.} Let $V_{i}=\{(v,i):v\in V(G)\}$ for $i=1,2,...,n$. Then $V(G\times C_{n})=\cup_{i= 1}^n {V_{i}}$. Let $X_{\frac{i+1}{2}}=\{(x,i):x\in X\}$ and $Y_{\frac{n+i+1}{2}}=\{(y,i):y\in Y\}$ for $i=1,3,...,n-1$. Let $X_{\frac{n+i}{2}}=\{(x,i):x\in X\}$ and $Y_{\frac{i}{2}}=\{(y,i):y\in Y\}$ for $i=2,4,...,n$. Then $V_{i}=X_{\frac{i+1}{2}}\cup Y_{\frac{n+i+1}{2}}$ for $i=1,3,...,n-1$ and $V_{i}=X_{\frac{n+i}{2}}\cup Y_{\frac{i}{2}}$ for $i=2,4,...,n$. Let
$H_{i}$ and $H_{i}^{'}$ be subgraphs of $G\times C_{n}$ induced by $X_{i}\cup Y_{i}$ and $X_{i+1}\cup Y_{i}$, respectively for $i=1,2,...,n$. For each $i$, $H_{i}$ and $H_{i}^{'}$ are isomorphic to $G$. Note that the two isomorphic components $G_{1}$ and $G_{2}$ are $\cup_{i= 1}^\frac{n}{2} {\left({{H_i\cup H_{i}^{'}}}\right)}$ and $\cup_{i= \frac{n}{2}+1}^n {\left({{H_i\cup H_{i}^{'}}}\right)}$.  Thus, by Lemma 2.2, both $G_{1}$ and $G_{2}$ (shown in Figure 3) are super-$\kappa$. $\Box$

\begin{figure}[htbp]
  \centering
  \includegraphics[width=13cm]{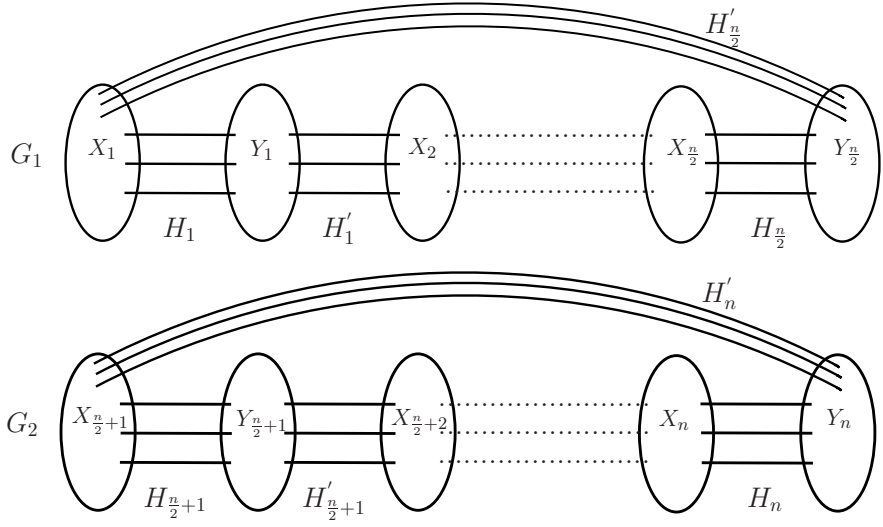}\\
  \caption{The graphs $G_{1}$ and $G_{2}$}
  \label{1}
\end{figure}

\begin{thm}
Let $G$ be a connected non-bipartite graph and $n\geq6$ be an even integer. If $\kappa(G\times K_{2})>\frac{4}{n}\delta(G)$, then $G\times C_{n}$ is super-$\kappa$.
\end{thm}

\noindent{\bf Proof.} Let $V_{i}=\{(v,i):v\in V(G)\}$ for $i=1,2,...,n$. Then $V(G\times C_{n})=\cup_{i= 1}^n {V_{i}}$. Let $H_{\frac{i+1}{2}}$ be the subgraph of $G\times C_{n}$ induced by $V_{i}\cup V_{i+1}$ for $i=1,3,...,n-1$. Let $H_{\frac{i}{2}}^{'}$ be the subgraph of $G\times C_{n}$ induced by $V_{i}\cup V_{i+1}$ for $i=2,4,...,n$. Then $H_{i}$ and $H_{i}^{'}$ are isomorphic to the bipartite graph $G\times K_{2}$ with bipartitions $(V_{2i-1}, V_{2i})$ and $(V_{2i}, V_{2i+1})$, respectively for $i=1,2,...,\frac{n}{2}$.
Note that $G\times C_{n}=\cup_{i= 1}^\frac{n}{2} {\left({{H_i\cup H_{i}^{'}}}\right)}$ (see Figure 4) and $\delta(G\times K_{2})=\delta(G)$.   By Lemma 2.2, $G\times C_{n}$ is super-$\kappa$. $\Box$

\begin{figure}[htbp]
  \centering
  \includegraphics[width=13cm]{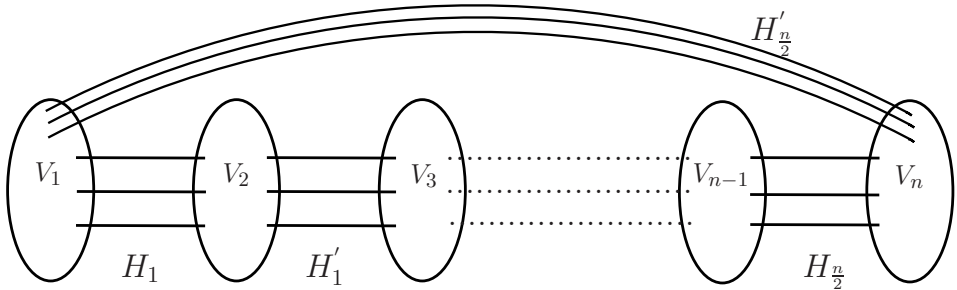}\\
  \caption{The graph $G\times C_{n}$ when $G$ is non-bipartite and $n$ is an even}
  \label{1}
\end{figure}

\begin{thm}
Let $G$ be a connected non-bipartite graph and $n\geq7$ be an odd integer. If $\kappa(G\times K_{2})>\frac{4}{n-1}\delta(G)$, then $G\times C_{n}$ is super-$\kappa$.
\end{thm}

\noindent{\bf Proof.} Let $V_{i}=\{(v,i):v\in V(G)\}$ for $i=1,2,...,n$. Then $V(G\times C_{n})=\cup_{i= 1}^n {V_{i}}$. Let $H_{\frac{i+1}{2}}$ be the subgraph of $G\times C_{n}$ induced by $V_{i}\cup V_{i+1}$ for $i=1,3,...,n$. Let $H_{\frac{i}{2}}^{'}$ be the subgraph of $G\times C_{n}$ induced by $V_{i}\cup V_{i+1}$ for $i=2,4,...,n-1$. Then $H_{i}$ is isomorphic to the bipartite graph $G\times K_{2}$ with bipartition  $(V_{2i-1}, V_{2i})$ for $i=1,2,...,\frac{n+1}{2}$ and $H_{i}^{'}$ is isomorphic to the bipartite graph $G\times K_{2}$ with bipartition  $(V_{2i}, V_{2i+1})$ for $i=1,2,...,\frac{n-1}{2}$.
Note that $G\times C_{n}=(\cup_{i= 1}^\frac{n-1}{2} {\left({{H_i\cup H_{i}^{'}}}\right)})\cup H_{\frac{n+1}{2}}$ (see Figure 5) and $\delta(G\times K_{2})=\delta(G)$. By similar arguments as the proof in Lemma 2.2, we can prove that $G\times C_{n}$ is super-$\kappa$. $\Box$

\begin{figure}[htbp]
  \centering
  \includegraphics[width=13cm]{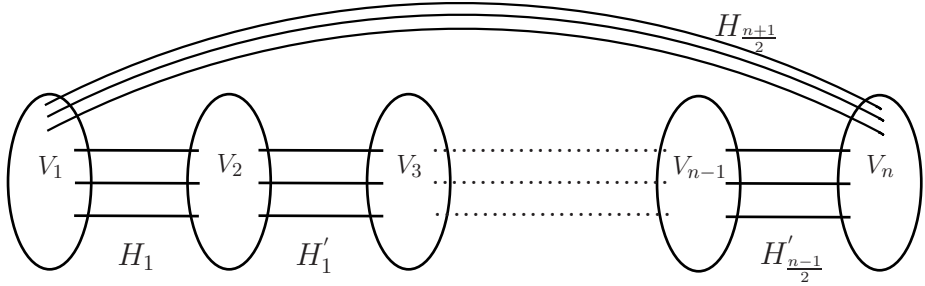}\\
  \caption{The graph $G\times C_{n}$ when $G$ is non-bipartite and $n$ is an odd}
  \label{1}
\end{figure}

\begin{thm} (\cite{Sonawane}) If $H$ is the direct product of $k\geq1$ odd cycles, then $\kappa(H\times K_{2})=2^{k}$.
\end{thm}


Combing Theorem 3.9 with Theorems 3.7 and 3.8 respectively, we have the following two corollaries.

\begin{cor} Let $G$ be the direct product of $k\geq1$ odd cycles and $n\geq6$ be an even interger. Then $G\times C_n$ is super-$\kappa$.
\end{cor}

\noindent{\bf Proof.} Let $G=C_{l_{1}}\times C_{l_{2}}\times\cdot\cdot\cdot\times C_{l_{k}}$. By Theorem 3.9, $\kappa(G\times K_{2})=2^{k}=\delta(G)$. Then $G\times C_{n}$ is super-$\kappa$ by Theorem 3.7. $\Box$

\begin{cor} Let $G$ be the direct product of $k\geq1$ odd cycles and $n\geq7$ be an odd interger. Then $G\times C_n$ is super-$\kappa$.
\end{cor}

\noindent{\bf Proof.} Let $G=C_{l_{1}}\times C_{l_{2}}\times\cdot\cdot\cdot\times C_{l_{k}}$. By Theorem 3.9, $\kappa(G\times K_{2})=2^{k}=\delta(G)$. Then, by Theorem 3.7, $G\times C_{n}$ is super-$\kappa$. $\Box$

\section{Concluding Remarks}

Motivated by the results on the connectivity of the direct product of graphs and cycles in \cite{Sonawane}, we focus on studying the super connectedness of the direct product of graphs and cycles in this paper. By using a key lemma we obtained in Section 2, we give some sufficient conditions for the direct product of a graph and a cycle to be super connected. However, there are few results on the connectivity of the direct product of two general graphs. So we are going to explore the connectivity and super connectedness of the direct product of two general graphs in the future.


\end{sloppypar}

\end{document}